\RequirePackage{fix-cm}
\documentclass[smallextended]{svjour3}
\smartqed
\usepackage[utf8]{inputenc}
\usepackage{tgpagella}
\usepackage{tcolorbox}
\usepackage{listings}
\usepackage{textcomp}
\usepackage{lineno}
\usepackage[section]{placeins}
\makeatletter
\AtBeginDocument{%
  \expandafter\renewcommand\expandafter\subsection\expandafter{%
    \expandafter\@fb@secFB\subsection
  }%
}
\makeatother

\usepackage{graphicx}
\usepackage[colorlinks]{hyperref}
\hypersetup{
     colorlinks   = true,
     citecolor    = blue
}
\usepackage{amsmath, amssymb, mathtools}
\usepackage{bm}
\usepackage{algorithmic}
\usepackage{algorithm}
\usepackage{geometry}
\usepackage{mathptmx}      
\usepackage{tabularx}

\usepackage{setspace}
\usepackage{xcolor}
\usepackage[sort,comma,numbers]{natbib}
\journalname{Computational Geosciences}

\begin{document}
\title{ Hybrid Gaussian-cubic radial basis functions for scattered data interpolation}
\author{Pankaj K Mishra  \and\\ Sankar K Nath \and\\ Mrinal K Sen \and \\Gregory E Fasshauer
}

\authorrunning{Mishra, et al.}

\institute{P. K. Mishra \at
              Former: Department of Geology and Geophysics, Indian Institute of Technology Kharagpur, India \\
              Current: Department of Mathematics, Hong Kong Baptist University, Kowloon Tong, Hong Kong.\\
              \email{pkmishra@gg.iitkgp.ernet.in}           
           \and
           S. K. Nath \at
           Department of Geology and Geophysics, Indian Institute of Technology Kharagpur, India \\
           \email{nath@gg.iitkgp.ernet.in}
           \and
           M. K. Sen \at
           Jackson School of Geosciences, University of Texas at Austin, USA\\
           \email{msentx@gmail.com}
           \and
           G. E. Fasshauer \at
           Department of Applied Mathematics and Statistics, Colorado School of Mines, USA\\
           \email{fasshauer@mines.edu }
}

\date{Received: 00-00-0000 / Accepted: 00-00-0000}
\maketitle
\def\linenumberfont{\normalfont\small\sffamily}
\begin{abstract}
\large
Scattered data interpolation schemes using kriging and radial basis functions (RBFs) have the advantage of being meshless and dimensional independent; however, for the data sets having insufficient observations, RBFs have the advantage over geostatistical methods as the latter requires variogram study and statistical expertise.
Moreover, RBFs can be used for scattered data interpolation with very good convergence, which makes them desirable for shape function interpolation in meshless methods for numerical solution of partial differential equations. For interpolation of large data sets, however, RBFs in their usual form, lead to solving an ill-conditioned system of equations, for which, a small error in the data can cause a significantly large error in the interpolated solution. In order to reduce this limitation, we propose a hybrid kernel by using the conventional Gaussian and a shape parameter independent cubic kernel. Global particle swarm optimization method has been used to analyze the optimal values of the shape parameter as well as the weight coefficients controlling the Gaussian and the cubic part in the hybridization. Through a series of numerical tests, we demonstrate that such hybridization stabilizes the interpolation scheme by yielding a far superior implementation compared to those obtained by using only the Gaussian or cubic kernels. The proposed kernel maintains the accuracy and stability at small shape parameter as well as relatively large degrees of freedom, which exhibit its potential for scattered data interpolation and intrigues its application in global as well as local meshless methods for numerical solution of PDEs.
\keywords{radial basis function \and multivariate interpolation \and particle swarm optimization \and spatial data analysis.}
\subclass{MSC 65 \and MSC 68 }
\end{abstract}

\section{Introduction}
\large
\noindent Generally used techniques for multivariate approximation were polynomial interpolation and piece-wise polynomial splines till 1971, Rolland Hardy proposed a new method using a variable kernel for each interpolation point, assembled as a function, depending only on the radial distance from the origin or any specific reference point termed as `\textit{center}', \citep{Hardy1971}. Such a kernel is known as radial kernel or radial basis functions (RBFs). In 1979, Richard Franke studied many available approaches for scattered data interpolation and established the utility of RBFs over many other schemes  {\color{blue}\citep{franke1979}}. Since RBF approximation methods do not require to be interpolated over tensor grids using rectangular domains, they have been proven to work effectively, where the polynomial approximation could not be precisely applied  {\color{blue}\citep{Sarra201168}}.

Global interpolation methods based on radial basis functions have been found to be efficient for surface fitting of scattered data sampled at $s$-dimensional scattered nodes. For a large number of samples, however, such interpolations lead to the solution of ill-conditioned system of equations  {\color{blue}\citep{Fass2009,Forn2011,Fass2012,Lin20122,Chen2014}}. This ill-conditioning occurs due to the global nature of the RBF interpolation, where the interpolated value at each node is influenced by all the node points in the domain providing full matrices that tend to become more ill-conditioned as the shape parameter gets smaller \citep{Driscoll2002}.

Several approaches have been proposed to deal with ill-conditioning in global interpolation using RBFs: Kansa and Hon performed a series of numerical tests using various schemes like replacement of global solvers by block partitioning or LU decomposition, use of matrix preconditioners, variable shape parameters based on the function curvature, multizone methods, and node adaptivity which minimizes the required number of nodes for the problem  {\color{blue}\citep{KansaHon2002}}. Some other approaches to deal with the ill-conditioning in the global RBF interpolation are: accelerated iterated approximate moving least squares  {\color{blue}\citep{Fass2009}}, random variable shape parameters  {\color{blue}\citep{Sarra20091239}}, Contour-Pad\'e  and RBF-QR algorithms  {\color{blue}\citep{Forn2011}}, series expansion of Gaussian RBF  {\color{blue}\citep{Fass2012}}, and regularized symmetric positive definite matrix factorization \citep{Sarra2014}. The Contour-Pad\'e approach is limited to few degrees of freedom only. According to Fornberg and Flyer  {\color{blue}\citep{Fornberg2013627,Forn2015}}, the current best (in terms of computation time) is the RBF-GA algorithm which is limited to the Gaussian RBFs only. RBF-GA is a variant of the RBF-QR technique developed by Fornberg and Piret  {\color{blue}\citep{Fornberg200760}}. There are other modern variants: the Hilbert-Schmidt SVD approach developed by Fasshauer and McCourt which can be shown to be equivalent to RBF-QR, but approaches the problem from the perspective of Mercer's theorem and eigenfunction expansions  {\color{blue}\citep{Fass2015}}. Another approach has been proposed by DeMarchi what he termed as the Weighted SVD method, which works with any RBF, but requires a quadrature/cubature rule, and only partially offsets ill-conditioning  {\color{blue}\citep{DeMarchi20131}}. Recently Kindelan et al. proposed an algorithm to study RBF interpolation with small parameters, based on Laurent series of the inverse of the RBF interpolation matrix {\color{blue} \citep{Gonzale2015,Kindelan2016}}. Another parallel approach to deal with the ill-conditioning problem is to use the local approximations like RBF-FD; however in this paper, we consider only the global approximations to focus on the effect of the proposed approach.

In this paper, we propose a hybrid radial basis function (HRBF) using a Gaussian and a cubic kernel, which significantly improves the condition of the system matrix avoiding the above mentioned ill-conditioning in RBF interpolation schemes. In order to examine the effect of various possible hybridization with these two kernels, we use the global particle swarm optimization for deciding the optimal weights for both kernels as well as the optimal value of the shape parameter of the Gaussian kernel. Through a number of numerical tests, we discuss the advantages of the proposed kernel over the Gaussian and the cubic kernels, when used individually. We also test the proposed kernel for the interpolation of a synthetic topographical data near a normal fault, at a relatively large number of desired locations, which could not be interpolated with only the Gaussian kernel, due to the ill-conditioning issue.

The rest of the paper is structured as follows. In section 2, we give a brief introduction of the radial basis functions and the fundamental interpolation problem. In section 3, after briefly arguing over the need of RBFs for interpolation problems and its advantage over geostatistical approaches for the same, we introduce our proposed hybrid kernel along with the motivation behind such a hybridization. Section 4 contains the description of polynomial augmentation in the hybrid kernel and additional constrains to ensure the non-singularity of the system. In order to get the optimal combination of the parameters introduced in section 3, we use a global particle swarm optimization algorithm, which is briefly explained in section 5. Through numerical tests, we check the linear reproduction property of the proposed hybrid kernel and application of this kernel for 2-D interpolation problems in sections 6.1 and 6.2 respectively. We analyze the eigenvalues of the interpolation matrices obtained by using the proposed hybrid kernel in section 6.3. In section 6.4, a comparative study of two different cost functions is presented. In order to bring in some complexity in the analysis, in section 6.5, the proposed hybrid radial basis function is used for interpolation of a 2-D geophysical data and its comparison with some other interpolation methods like linear, cubic, ordinary and universal kriging and Gaussian radial basis interpolation. Finally, we discuss the computational cost of the present method, followed by our conclusions.

\section{Radial basis functions and spatial interpolation}
RBF and geostatistical approaches for scattered data interpolation are quite similar and provide almost equally good results, in general; however, for the data set which has very small number of measurements at locations making the spatial correlation very difficult, the application of geostatistical methods offers several challenges.
In such cases, variographic study is quite difficult to perform, which is a primary requirement in understanding the data for proper application of statistical tools. In 2009, Cristian and Virginia Rusu performed several interpolation tests with real data sets and discussed some of the advantages of using RBFs for scattered data interpolation  {\color{blue}\citep{Rusu2006}}. For completeness, we briefly define the RBF and the fundamental interpolation problem.

\noindent \textbf{Definition 2.1.} A function $\bm{\Phi}: \mathbb{R}^s \rightarrow \mathbb{R}$ is said to be radial if there exists a univariate function $\phi : [0, \infty) \rightarrow \mathbb{R}$ such that
\begin{equation}
\bm{\Phi} \left(\bm{x} \right) = \phi(r),   \qquad r=\parallel\bm{x}\parallel.
\end{equation}

\noindent$\parallel\cdot\parallel $ here, represents Euclidean norm. Some commonly used RBFs have been listed in Table \ref{tab:rbflist}. The constant $\epsilon$ is termed as the shape parameter of the corresponding RBF.
\begin{table*}
  \large
  \begin{tabular*}{250pt}{ll}
    \hline
    RBF Name & Mathematical Expression \\
    \hline
    Multiquadric              & $\phi(r)= (1+(\epsilon r)^2)^{1/2}$ \\
    Inverse multiquadric      & $\phi(r)= (1+(\epsilon r)^2)^{-1/2} $\\
    Gaussian                    & $\phi(r)=e^{-(\epsilon r)^2}$\\
    Thin plate spline            & $\phi(r)=r^2 \log(r) $\\
    Cubic                        & $\phi(r)=r^3$\\
    Wendland's                  & $\phi(r)=(1-\epsilon r)^{4}_{+}(4\epsilon r+1)$\\
    \hline
\end{tabular*}
 \caption{Typical RBFs and their expressions}
   \label{tab:rbflist}
\end{table*}
\noindent For completeness, we briefly define a general interpolation problem here. Let us assume that we have some measurements $y_j\in \mathbb{R}$ at some scattered locations $\bm{x}_j \in \mathbb{R}^s$. For most of the benchmark tests, we assume that these measurements $\left(y_j\right)$ are obtained by sampling some test function $\bm{f}$ at locations $\bm{x}_j$. \\
\noindent \textbf{Problem 2.1.} Given the data $y_j \in \mathbb{R}$ at locations $\bm{x}_j \in \mathbb{R}^s$ $(j=1,2,..N)$, find a continuous function $\mathcal{F}$ such that,
\[\mathcal{F}(\bm{x}_j) = y_j, \qquad  j=1,2,...N.\]
\noindent Given the set of scattered centers $\bm{x}_j$, an RBF interpolant can be written as
\begin{eqnarray}{\label{eq:one}}
\mathcal{F}(\bm{x}) = \sum_{j=1}^{N} c_j \phi (\parallel \bm{x}-\bm{x_j}\parallel),
\end{eqnarray}
where $\phi (\parallel \bm{x}-\bm{x_j}\parallel)$ is the value of the radial kernel, $\parallel \bm{x}-\bm{x_j}\parallel$ is the Euclidean distance between the observational point and the center, and $c_j$ are the unknown coefficients which are determined by solving a linear system of equations depending on the interpolation conditions. The system of linear equations for above representation can be written as \\
\begin{eqnarray}
\begin{bmatrix}
 \phi (\parallel \bm{x}_1-\bm{x}_1\parallel)          &  \cdots     & \phi (\parallel \bm{x}_1-\bm{x}_N\parallel)   \\
\phi (\parallel \bm{x}_2-\bm{x}_1\parallel)         & \cdots & \phi (\parallel \bm{x_2}-\bm{x}_N\parallel) \\
 \vdots & \vdots & \ddots &         \vdots  \\
 \vdots & \vdots &        &  \vdots \\
\phi (\parallel \bm{x_N}-\bm{x}_1\parallel)        & \cdots &\phi (\parallel \bm{x}_N-\bm{x}_N\parallel)
\end{bmatrix}
\begin{bmatrix}
c_1 \\
c_2\\
\vdots\\
\vdots \\
c_N
\end{bmatrix} =
\begin{bmatrix}
f(\bm{x_1}) \\
f(\bm{x_2})\\
\vdots\\
\vdots \\
f(\bm{x_N})
\end{bmatrix}.
\end{eqnarray}

\section{Hybrid Gaussian-cubic kernels}
\noindent The stability and accuracy of the RBF interpolation depends on certain aspects of the algorithm, and the data involved. For example, if the scattered data comes from a sufficiently smooth function, the accuracy of the interpolation using the Gaussian RBF will increase relatively more as we increase the number of data points used, than those obtained with other kernels. Use of positive definite kernels like the Gaussian radial basis function assures the uniqueness in the interpolation. This means that the system matrix of interpolation is non-singular even if the input data points are very few and poorly distributed making the Gaussian RBF a popular choice for interpolation and numerical solution of PDEs. Since the accuracy and stability of Gaussian radial basis function interpolation mostly depends on the solution of the system of linear equations, the condition number of the system matrix plays an important role. It has been seen that the system of linear equations is severely ill-conditioned either for large number of data points in the domain or for very small shape parameter i.e., flat radial basis. A ``small" shape parameter becomes typical for sparse data. Also, to deal with coarse sampling and a low-level information content, the RBFs must be relatively flat. The Gaussian radial basis function leads to ill-conditioned system when the shape parameter is small.

Cubic radial basis function $(\phi (r) = r^3)$, on the other hand, is an example of finitely smooth radial basis functions. Unlike the Gaussian RBF, it is free of shape parameter which excludes the possibility of ill-conditioning due to small shape parameters. However, using a cubic RBF only in the interpolation might be problematic since the resulting linear system may become singular for certain point locations. Another difference between cubic RBFs and Gaussians lies in their approximation power. This means that assuming the scattered data comes from a sufficiently smooth function, for a Gaussian RBF the accuracy of an interpolant will increase much more rapidly than the cubic RBF, as the number of used data points are increased. On the other hand, cubic RBFs can provide more stable and better converging interpolations for some specific data set but the risk of singularity will always be there associated with a typical node arrangement depending on the data type.

From what we have discussed above, it is certain that both the Gaussian and cubic RBFs have their own advantage for scattered data interpolation. In order to make the interpolation more flexible, we propose a hybrid basis function using a combination of both the Gaussian and the cubic radial basis function as given by
\begin{eqnarray}
\phi(r) = \alpha e^{ -(\epsilon r)^2} +\beta r^3.
\end{eqnarray}

\begin{figure}[hbtp]
\centering
\includegraphics[scale=0.8]{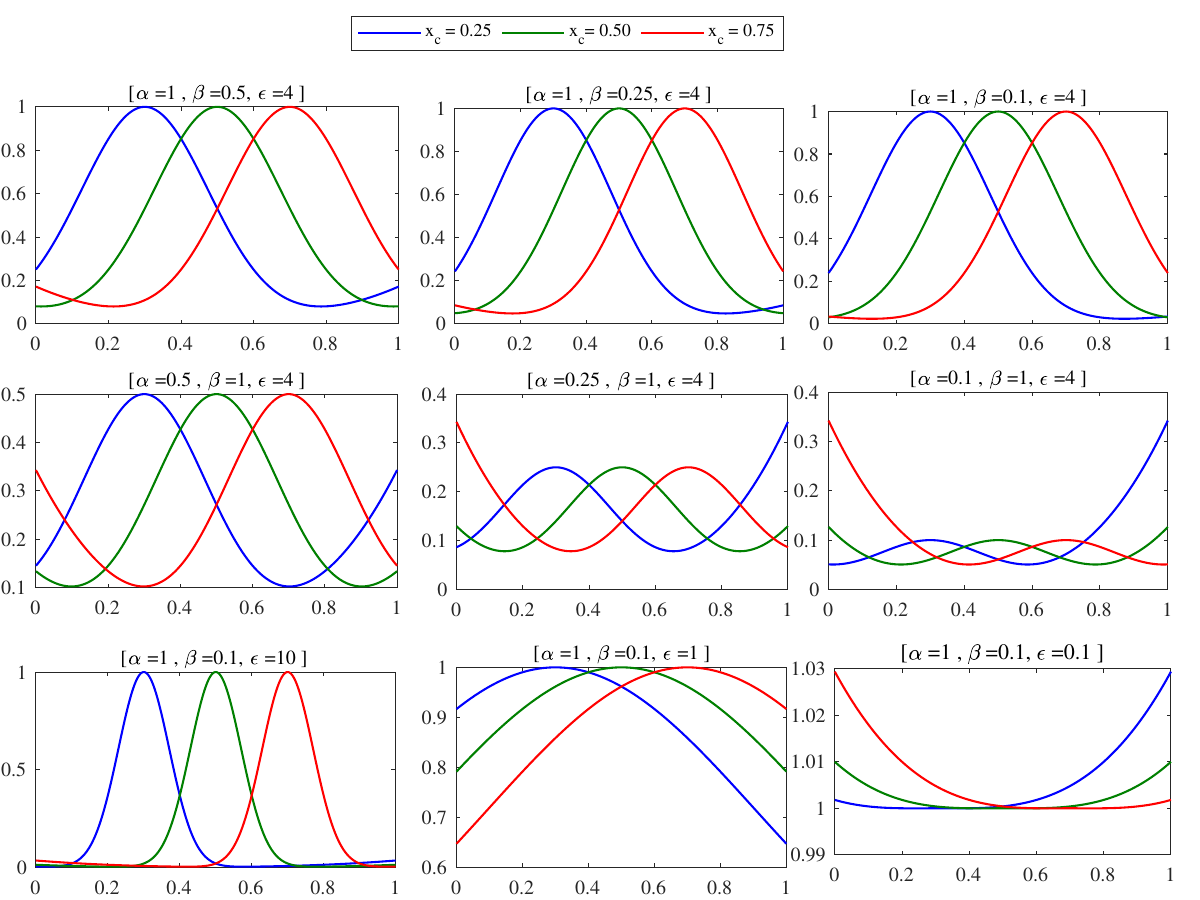}
\caption{\large 1D Plots of the hybrid kernels with different combinations of the parameters (x-axis: spatial location). Kernels have been plotted with three different centers ($x_c$).}
\label{1dkernel}
\end{figure}

\noindent Here, $\epsilon$ is the usual shape parameter associated with the RBFs. We have introduced two weights $\alpha$ and $\beta$, which control the contribution of the Gaussian and the cubic part in the hybrid kernel, depending upon the type of problem and the input data points, to ensure the optimum accuracy and stability. Such a combination would seem appealing for different reasons, such as (1) the involvement of the cubic kernel in the hybridization helps conditioning when using low values of shape parameter, which likely is good for accuracy and (2) supplementary polynomials (of low order) can help cubics at boundaries and also improve (derivative) approximations of very smooth functions (such as constants). Figure \ref{1dkernel} contains the 1-D plots of the proposed kernel with some arbitrarily chosen parameters.  In order to find the best combination of $\epsilon$, $\alpha$, and $\beta$, corresponding to maximum accuracy, we use the global particle swarm optimization algorithm.

\section{Polynomial augmentation and linear reproduction}
Application of radial basis functions for multivariate data interpolation is strongly supported by the limitation of polynomials for the same. In general, the well-posedness of multivariate polynomial interpolation can only be assured if the data sites are in certain special locations. However, polynomial augmentation in radial basis function has been reported to improve the convergence. A desirable property of a typical interpolant is linear reproduction, \emph{i.e.}, if the data is sampled from a linear function, the interpolation should reproduce it exactly. The interpolation using the Gaussian radial basis function does not reproduce the simple linear polynomials. Such polynomial reproduction is recommended especially when the interpolation is intended to be used in numerical solution of partial differential equations. The ability to reproduce the linear function is termed as the ``\emph{patch test}''. The Gaussian, and other radial basis functions, which are integrable, do not reproduce linear polynomials  {\color{blue}\citep{Sback1994,Fassbook2007}}. However, a polynomial augmentation in the Gaussian radial basis function reproduces the linear polynomial to the machine precision  {\color{blue}\cite{Fassbook2007}}. Fornberg and Flyer have also shown the improvement due to polynomial augmentation in cubic and quintic radial basis function  {\color{blue}\cite{Fornberg2002,Barnett2015}}.

It is well known that the Gaussian kernel is positive definite. This implies that the Gaussian kernel is conditionally positive definite of any order too. Any non-negative linear combination of conditionally positive definite kernels (of the same order) is again conditionally positive definite (of the same order). Therefore, adding a linear polynomial to the hybrid kernel expansion ensures that the resulting augmented interpolation matrix is non-singular; just as it does for the cubic kernel.

In order to reproduce linear functions, we add low order polynomial function, \emph{i.e.}, $\bm{x}\mapsto 1$,$\bm{x}\mapsto x$, and $\bm{x}\mapsto y$ to the proposed kernel. By doing so, we have added three more unknowns to the interpolation problem, which makes the total number of unknowns equal to $N+3$. The polynomial augmentation, explained here, is for the 2D case only. Linear polynomial augmentation in other dimensions works analogously. The approximation now becomes,
\begin{eqnarray}
\mathcal{F}(\bm{x})  = \sum_{k=1}^{N} c_k \phi(\parallel \bm{x}-\bm{x_k}\parallel)+c_{N+1}+c_{N+2}x+c_{N+3}y \qquad \bm{x}=(x,y)\in \mathbb{R}^2,
\end{eqnarray}
At this point, we have $N$ interpolation conditions, i.e.,
\begin{eqnarray}
\mathcal{F}(\bm{x}_j) = \bm{f}(\bm{x}_j) \qquad  j=1,2,...N.
\end{eqnarray}
In accordance with \cite{Fassbook2007,Sback1994,Fornberg2002}, we add three additional conditions as follows:
\begin{eqnarray}
\sum_{k=1}^N c_k =0, \qquad \sum_{k=1}^N c_k x_k =0, \qquad \sum_{k=1}^N c_k y_k =0.
\end{eqnarray}
These additional constraint lead to a square (non-singular) system of linear equations. Now we need to solve a system of the form
\begin{eqnarray}
\begin{bmatrix} \mathbf{A} & \mathbf{P} \\ \mathbf{P}^T & O \end{bmatrix}
\begin{bmatrix} \mathbf{c} \\ \mathbf{d}\end{bmatrix}
=
\begin{bmatrix} \mathbf{y}\\ \mathbf{0}\end{bmatrix},
\end{eqnarray}
where
\[\mathbf{A}_{j,k}= \phi(\parallel \bm{x}_j-\bm{x}_k\parallel),  \qquad j,k=1,...,N, \]
\[\mathbf{c}=[c_1,...,c_N]^T,\]
\[\mathbf{y}=[f(\bm{x}_1)....,f(\bm{x}_N)]^T,\]
\[\mathbf{P}_{j,k} =p_k(\bm{x}_j), \qquad p =[1\qquad x \qquad y], \]
$ \mathbf{0}$ is a zero vector of length 3, and $O$ is a zero matrix of size $3\times 3$.

\section{Parameter Optimization}
\noindent  The selection of a good shape parameter in an RBF-based algorithm has always been a prime concern. In the context of dealing with ill-conditioned systems in RBF interpolation and application, the work done by Fornberg and his colleagues established that there exists an optimal value of the shape parameter, which corresponds to the optimal accuracy, and that this value can be determined in a stable manner. The performance of the hybrid kernel varies with different combinations of parameters: $\alpha$, $\beta$, and $\epsilon$. Therefore, for a typical problem --- there must be an optimum combination of these parameters for which the performance of the hybrid kernel is the ``best''.  Since the hybrid kernel involves more than one parameter, we need a more general approach to find optimal kernel parameters. We propose an algorithm for parameter selection of the hybrid kernel by using the global particle swarm optimization (PSO) approach, which is discussed in detail in \textbf{Appendix B}.  We have considered the following two approaches to compute the objective function in the optimization, and thus define what we mean by ``best'':
\begin{description}
\item[$\bullet$] \textbf{RMS Error:} The root mean square error over $M$ evaluation points is computed according to the formula given by,
\begin{eqnarray}
\label{error}
E_{rms} = \sqrt{\frac{1}{M}\sum_{j=1}^{M} \left[\mathcal{F}(\xi_j) - f(\xi_j) \right]^2}
\end{eqnarray}
Where $\xi_j, (j=1,...,M)$ are the evaluation points and $E_{rms}$ is the error function which is to be optimized for the minimum values for a set of $\epsilon$, $\alpha$, and $\beta$.

\item[$\bullet$] \textbf{Leave-one-out-cross-validation:} In practical problems of scattered data interpolation, the exact solution is most likely to be unknown. In such situations, it is not possible to calculate the exact RMS error. Cross-validation is a statistical approach. The basic idea behind cross-validation is to separate the available data into two or more parts and test the accuracy of the involved algorithm. Leave-one-out-cross-validation (LOOCV) in particular, is a special case of an exhaustive cross-validation, \textit{i.e.}, leave-$p$-out-cross-validation with $p=1$. The LOOCV algorithm separates one datum from the whole dataset of size $N$ and calculates the error of the involved method using the rest of the data. This is done repeatedly, each time excluding different datum and the corresponding errors are stored in an array. Unlike the regular LOOCV approach which performs relative error analysis \cite{Who90}, we consider this LOOCV error as the objective function in the particle swarm optimization algorithm. A detailed explanation of LOOCV can be found in \citep{friedman2001elements,GETOOR2007,Rippa}, however, we briefly explain the LOOCV used in this paper in \textbf{Appendix A}.
\end{description}
\noindent The optimization problem here, can be written in the mathematical form as following:

\begin{tcolorbox}
\begin{equation*}
\begin{aligned}
& \underset{\alpha,\beta,\epsilon}{\text{minimize}}
& & O_f(\alpha,\beta,\epsilon) \\
& \text{subject to}
& & \epsilon \geq 0,\\
& & &  0 \leq \alpha \leq 1,\\
& &  &  0 \leq \beta \leq 1,
\end{aligned}
\end{equation*}
\begin{description}
\item where
\item[$\bullet$] $O_f$ is the objective function, which is a cost vector computed either through RMS error or LOOCV.
\item[$\bullet$] $\alpha, \beta$ and $\epsilon$ are the kernel parameters.
\end{description}
\end{tcolorbox}

\section{Numerical Tests}
\subsection{Linear Reproduction}
\begin{figure}[hbtp]
\centering
\includegraphics[scale=0.5]{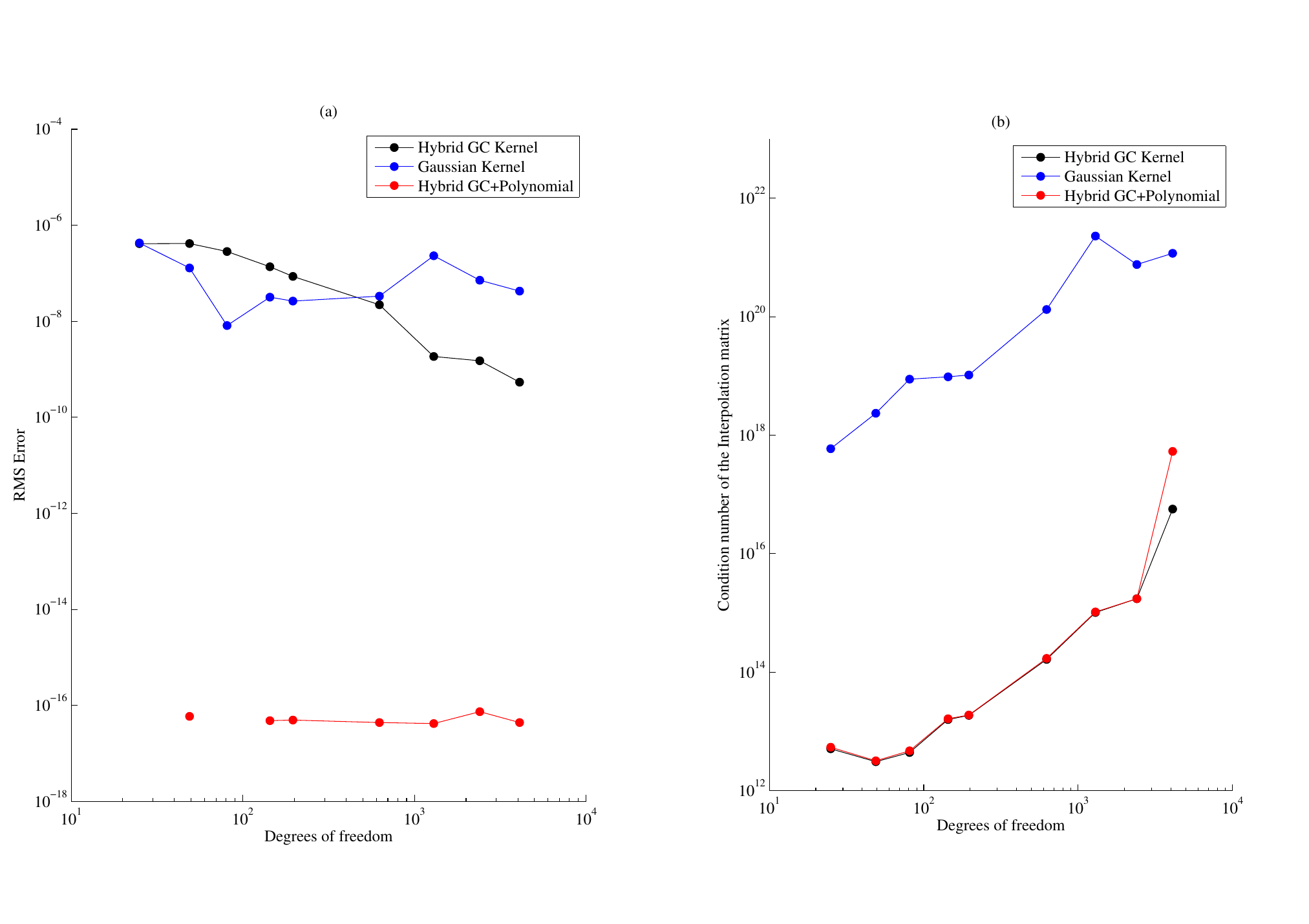}
\caption{Numerical test to check the linear reproduction property of the proposed hybrid Gaussian-cubic kernel. The RMS error convergence and the condition number variation with increasing degrees of freedom using the Gaussian, the hybrid and the hybrid kernel with polynomial augmentation have been compared. }
\label{fig:LRtest2}
\end{figure}

\noindent In this test, we examine the linear reproduction of the hybrid kernel and compare it with that of the Gaussian kernel. We consider the following linear function to sample the test data,
\[f(\mathbf{x})= \frac{x+y}{2}.\]
\noindent Table \ref{tab:LRtest} contains the results of particle swarm optimization for the aforementioned linear polynomial reproduction. The objective function, which is frequently referred as the ``cost function" in the RBF literature, is RMS error. In order to compare the convergence of linear reproduction test using the Gaussian kernel, the hybrid kernel and the hybrid kernel with polynomial augmentation, parameter optimization has been performed for various degrees of freedoms, i.e., $[ 25, 49, 144, 196,  625, 1296, 2401, 4096]$. The RMS errors for various degrees of freedom using the hybrid kernel, the Gaussian kernel and the hybrid kernel with polynomial augmentation has been denoted as $E_{H}$, $E_{G}$, and $E_{H+P}$ respectively. Although we have used PSO for the optimization of the parameters here, it may not play a relevant role here for the fact that there are many (near) optimal results which reproduce the linear polynomial. It is well known that the Gaussian kernel alone does not reproduce any order of polynomial and needs polynomial augmentation for the same. This limitation is continued with the proposed hybrid kernel too. However, if we notice the convergence comparison for this test in Figure \ref{fig:LRtest2}(a), the proposed hybrid kernel offers better convergence as compared to the Gaussian kernel, which actually shows no convergence. The reason for such better convergence is the reduced condition number due to the proposed hybridization which could be seen in Figure \ref{fig:LRtest2}(b). Although the hybrid kernel does not reproduce the linear polynomial exactly, the condition number of the interpolation matrix is almost similar to those obtained by hybrid kernel with polynomial augmentation. Interestingly, a very small doping ($\beta$ values in Table 2) of the cubic kernel into the Gaussian kernel reduces the condition number of the interpolation, which is a solution to the ill-conditioned problem of the Gaussian kernel at large degrees of freedom.
\begin{table*}
\small
  \centering
  \begin{tabular*}{\textwidth}{l@{\extracolsep\fill}cccccc}
   \hline
    N & $\epsilon $ & $\alpha $ & $\beta$ & $E_{H}$ & $E_{G}$ & $E_{H+P}$ \\
    \hline
  25	& 0.1600	&0.9592	&$1.73e-09$	&$4.14e-07$     &$4.27e-07$	   	&$0.00e-00$	 \\

  49	&0.2151	&0.9216	&$2.03e-08$	&$4.15e-07$     &$1.29e-07$		&$6.00e-17$	\\

  81	&0.5580	&0.5790	&$3.29e-08$	&$2.83e-07$		&$8.14e-09$		&$0.00e-00$	\\

 144	&0.8633	&0.6994	&$4.58e-08$	&$1.36e-07$		&$3.18e-08$		&$4.81e-17$	\\

 196	&0.5911	&0.9277	&$1.34e-07$	&$8.57e-08$		&$2.63e-08$		&$4.96e-17$	\\

 625	&1.5290	&0.4350	&$8.89e-08$	&$2.20e-08$		&$3.33e-08$		&$4.40e-17$	\\

1296	&0.9397	&0.2791	&$8.46e-08$	&$1.85e-09$		&$2.31e-07$		&$4.18e-17$	\\

2401	&0.9450	&0.7403	&$6.27e-07$	&$1.50e-09$		&$7.14e-08$		&$7.43e-17$	\\

4096	&1.1183	&0.7590	&$2.34e-07$	&$5.37e-10$		&$4.24e-08$		&$4.40e-17$	\\
    \hline
\end{tabular*}
\caption{Linear reproduction with parameter optimization, using hybrid kernel with polynomial augmentation.}
\label{tab:LRtest}
\end{table*}
Another interesting observation in this numerical test is the small values of the shape parameter. This observation exhibits the efficiency of the proposed hybridization at small shape parameters, which explains the higher accuracy of linear reproduction. The flat kernels made possible by the hybridization are closer to polynomials than what is possible with the Gaussian kernel alone due to the more severe ill-conditioning.

\subsection{Franke's test}
\noindent Here, we perform a 2D interpolation test using the proposed hybrid kernel and its comparative study with the Gaussian and cubic spline. A benchmark test \citep{franke1979} for 2D interpolation is to sample and reconstruct Franke's test function which is given by,
\begin{eqnarray}
\label{Franke}
f(x,y) = f_1+f_2+f_3-f_4,
\end{eqnarray}
where
 \[f_1=0.75 e^{\left( -\frac{1}{4} \left( (9x-2)^2 + (9y-2)^2 \right)\right)},\]
\[ f_2=0.75 e^{\left( -\frac{1}{49} (9x+1)^2 + \frac{1}{4}(9y+1)^2 \right)}, \]
\[f_3=0.50 e^{\left( -\frac{1}{4} \left( (9x-7)^2 + (9y-3)^2 \right)\right)},\]
\[f_4=0.20 e^{\left( (9x-4)^2 - (9y-7)^2 \right)}.\]

\begin{figure}[hbtp]
\centering
\includegraphics[scale=0.4]{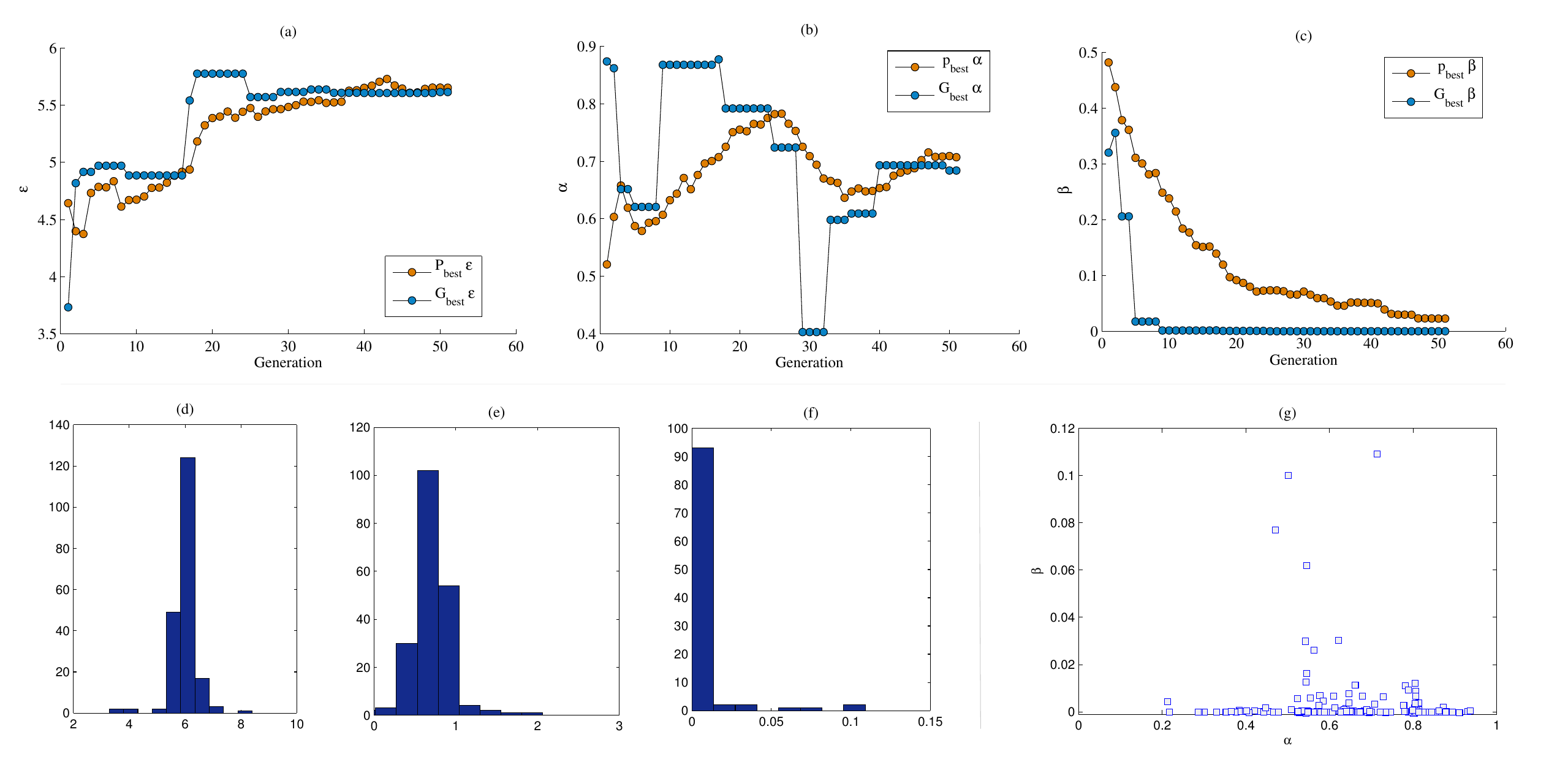}
\caption{(a)-(c) Convergence of $pbest$ and $gbest$ over generations for optimization of $\epsilon$, $\alpha$, and $\beta$ respectively. (d)-(f) frequency histogram of the optimized solution for 40 swarms over 5 iterations for $\epsilon$, $\alpha$, and $\beta$ respectively. (g) shows the
interdependence of $\alpha$ and $\beta$ for fixed value of shape parameter $\epsilon = 5.5434$.}
\label{fig:PSOconv}
\end{figure}

\begin{figure}[hbtp]
\centering
\includegraphics[scale=0.49]{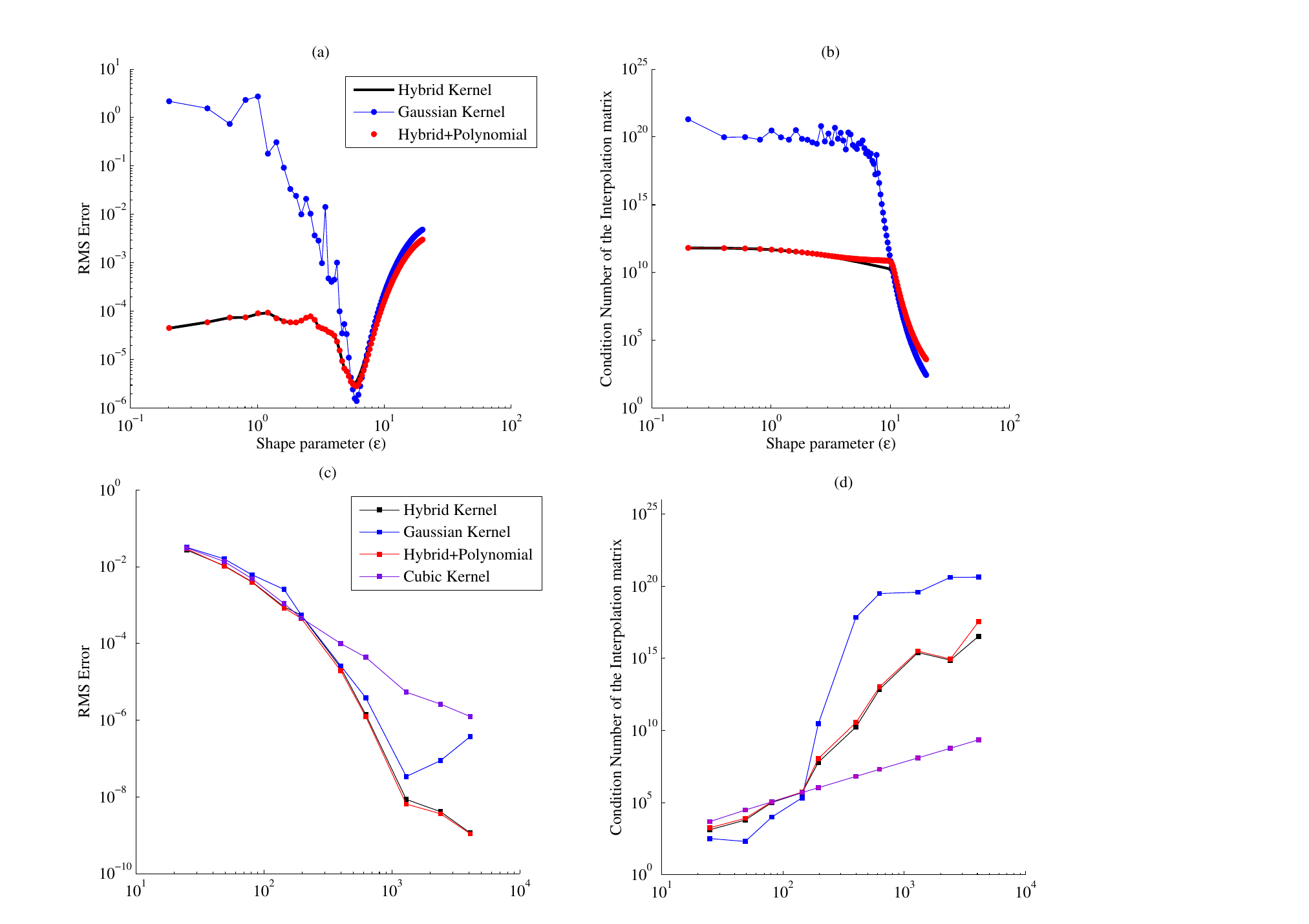}
\caption{Interpolation test with Franke's Test function using RMS Error as the objective function; (a)-(b) RMS error variation with the shape parameter for various kernel and (c)-(d) RMS error convergence for various kernels.}
\label{fig:Franke2D}
\end{figure}

\begin{figure}
\centering
\includegraphics[scale=0.4]{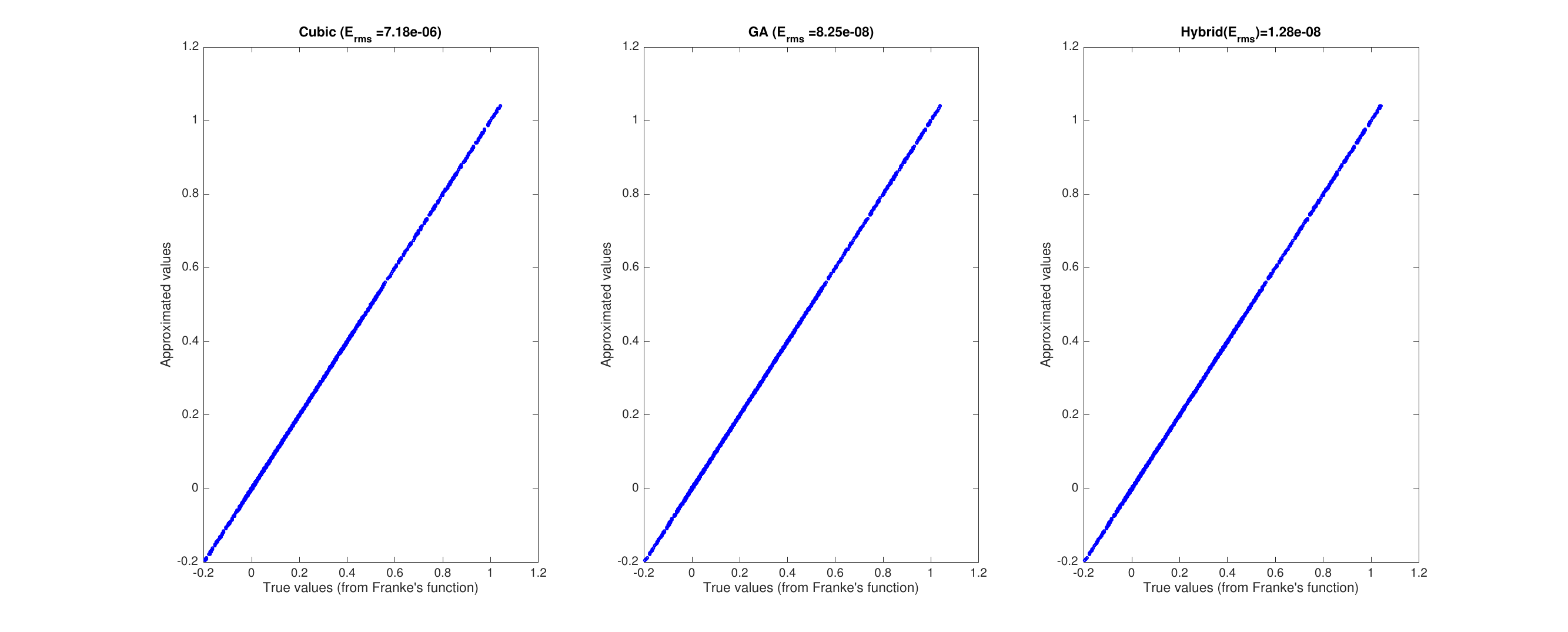}
\caption{Crossplot between true solution obtained from Franke's function and the approximated solution with different kernels with optimal parameters. (N=1225)}
\label{fig:crossplot}
\end{figure}

 The RMS error was kept as the objective function in the particle swarm optimization. We sample the Franke's test function at various number of data points, \textit{i.e.}, $ [ 25, 49, 144, 196, 625, 1296, 2401, 4096]$ and try to reconstruct it while finding the optimal parameter combination for the hybrid kernel. Figure \ref{fig:PSOconv} shows a typical particle swarm optimization procedure for this test when we try to reconstruct the Franke's test function using $625$ data points. The convergence of the $pbest$ and $gbest$ values of $\epsilon$, $\alpha$, and $\beta$ over the generations, are shown in Figures \ref{fig:PSOconv}(a)-(c). The swarm size for this test was $40$ and the optimization has been performed for 5 iterations. Figures \ref{fig:PSOconv}(d)-(f) show the histogram of the parameters' values acquired by each swarm over 5 iterations ($40\times 5=200$).
Table \ref{tab:PSO2DFranke} contains the results of this numerical test. We observe that the hybrid kernel performs very similar with or without the polynomial augmentation. Figure \ref{fig:Franke2D}(a) shows the RMS error variation with different shape parameters in this test with $625$ data points using hybrid kernel, hybrid kernel with polynomial augmentation, and only the Gaussian kernel. There are two interesting observations in Figure \ref{fig:Franke2D}(a). First, the optimal value of the shape parameter $\epsilon$ is almost the same for all the three kernels, which implies that the hybridization or even an additional polynomial augmentation does not affect the original optimal shape of the Gaussian part in the hybrid kernel. Also, Figure \ref{fig:Franke2D}(d) shows the interdependence of $\alpha$ and $\beta$ for which the corresponding $\epsilon$ does not change significantly. The observation here is that the optimization algorithm suggests various values of $\alpha$ in the solution implying that the weight of the Gaussian kernel could have any values between $0.2-0.9$. The second interesting observation in Figure \ref{fig:Franke2D}(a) is the performance of the hybrid kernel for very small values of shape parameters. The hybrid kernel with or without the polynomial augmentation is quite stable in this range too. The corresponding condition number variations are shown in Figure \ref{fig:Franke2D}(b). Like the previous numerical test of linear reproduction, the condition number here is also significantly reduced for the hybrid kernel, with or without polynomial augmentation, which is not surprising since the condition number of the matrix does not depend on the data values, only on the kernel chosen and the data locations. Figure \ref{fig:Franke2D}(c) shows the RMS error convergence with increasing degrees of freedom for the aforementioned kernels and an additional cubic kernel. The convergence of hybrid kernel with or without polynomial augmentation is better than that for the Gaussian kernel and far better than that of cubic kernel. The cubic kernel, however, has the least values of the condition numbers as shown in Figure \ref{fig:Franke2D}(d). The overall observation in this numerical test is that the hybrid kernel works better than the Gaussian and the cubic kernel as far as the RMS error convergence is concerned (Figure \ref{fig:crossplot}). Also, the polynomial augmentation is not required in the hybrid kernel unless the data comes from a linear function.
\begin{table*}
  \centering \small
  \label{test_kernelp2}
  \begin{tabular*}{\textwidth}{c@{\extracolsep\fill}cccccccc}
   \hline
    N & $\epsilon $ & $\alpha $ & $\beta$ & $E_{rms}$ & $\epsilon^p$ & $\alpha^p$ & $\beta^p$ & $E^{p}_{rms}$ \\
    \hline
    25	& 2.9432 &	$3.161e-01$     &   $4.661e-01$ & 	$2.724e-02$    & $3.7378$ & $8.253e-01$	& $ 2.544e-06$ &	$ 2.552e-02$\\
    49	& 4.8600 &    $1.138e-01	$   &  $ 8.603e-01$	& $ 1.070e-02$     & $5.0242$ & $1.633e-01$  & $	5.817e-01$ & $ 1.029e-02$\\
    81 & 5.1345 &  	$4.462e-02$	   &$ 9.316e-01$	    & $ 4.044e-03$ & $5.2688$ &$4.532e-02$   & $9.582e-01$ & $ 3.900e-03 $ \\
   144 & 6.2931 & 	$1.700e-02$ 	   &$ 8.494e-01$    & $	9.054e-04$  & $6.6797$ &$1.500e-02$   &$9.970e-01$ & $ 8.287e-04 $\\
   196 & 5.5800 &	$7.087e-02	$   & $9.445e-01$    &	$1.658e-04$     &$5.3149$  & $4.134e-01$ & $4.094e-01$ & $ 2.912e-04 $ \\
   400 & 5.5683 &	$4.500e-01$     &$ 4.649e-05$ 	&  $2.311e-05$ & $5.7856$     & $6.531e-01$ & $4.275e-07$ & $ 1.718e-05$\\
   625 & 5.5434 &	$6.749e-01$     & $4.915e-07 $   &	$1.400e-06$ & $5.7530$     & $7.584e-01$ & $1.342e-06$ & $ 1.716e-05$\\
  1296 & 6.2474 &    $7.880e-01$	 & $9.109e-09$	& $ 8.582e-09$ & $5.9265$     & $8.790e-01$ & $1.832e-09$ & $ 5.073e-09 $\\
  2401 & 6.0249 &	$5.600e-01$  & $2.503e-08$ &	$2.106e-09$ & $6.3070$        & $9.520e-01$ &$5.704e-09$ & $9.006e-10 $\\
  4096 & 5.7700&    $9.107e-01$ &  $7.090e-08	 $ & $1.150e-09$ & $5.9397$          &$6.548e-01$ & $1.756e-08$ & $ 7.730e-10$\\
    \hline
\end{tabular*}
\caption{Results of the parameter optimization test for 2-D interpolation using hybrid kernel. The superscript $p$ means that the parameter has been optimized with linear polynomial augmentation in the hybrid kernel.}
\label{tab:PSO2DFranke}
\end{table*}

\newpage
\subsection{Eigenvalue Analysis of Interpolation Matrices}
\begin{figure}[hbtp]
\centering
\includegraphics[scale=0.38]{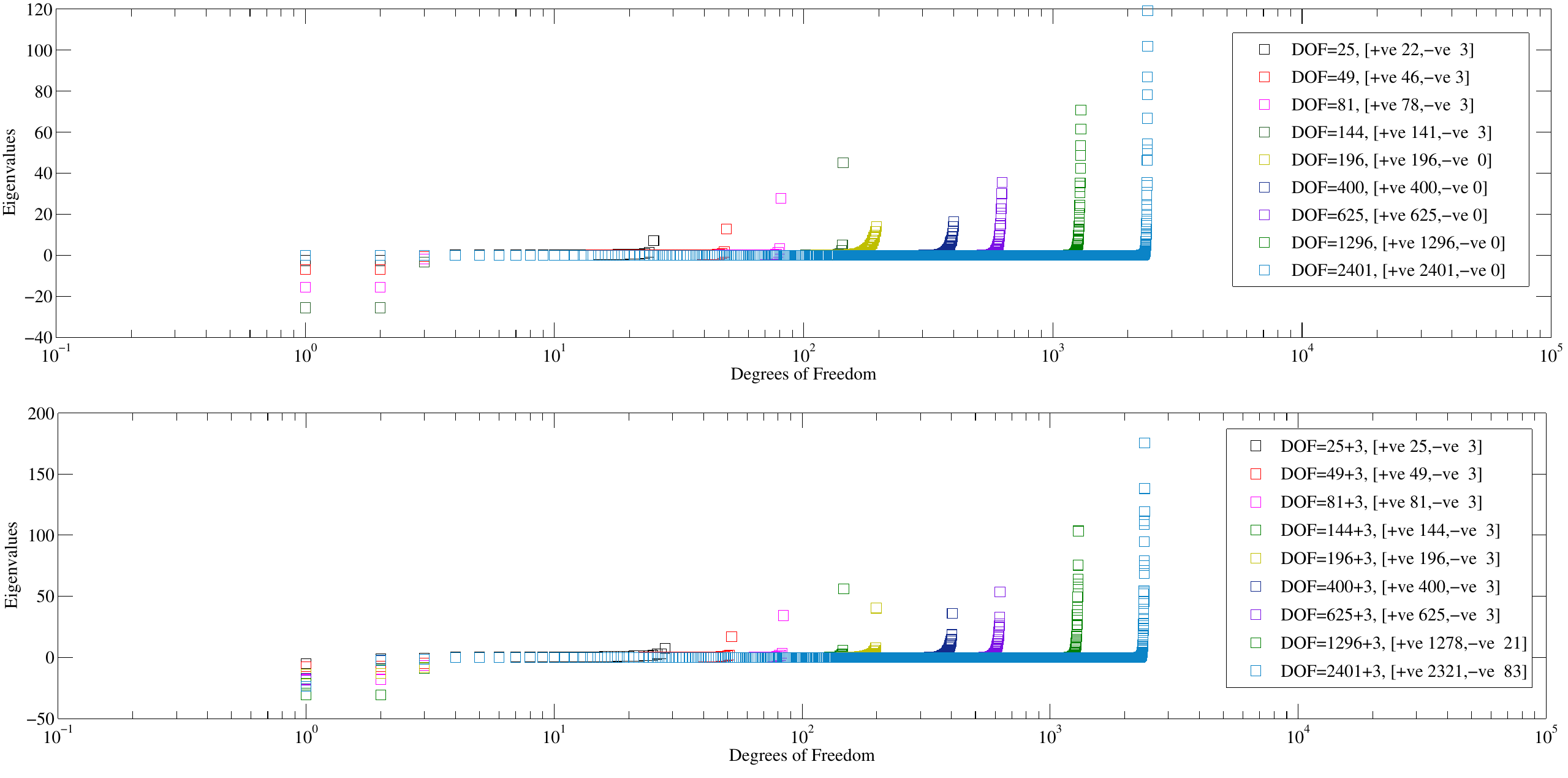}
\caption{Eigenvalue analysis of the interpolation matrix in 2-D interpolation using Franke's test function and corresponding optimal values of $\epsilon$, $\alpha$, and $\beta$. The objective function for the optimization was RMS Error. The eigenvalue spectra of the interpolation matrix has been shown for the hybrid kernel (top) and the hybrid kernel with polynomial augmentation (bottom).}
\label{fig:eigentest}
\end{figure}

\noindent In this numerical test, we visualize the eigenvalue spectra of the interpolation matrices for the hybrid kernel with and without polynomial augmentation. Figure \ref{fig:eigentest} shows the eigenvalue analysis of the interpolation matrix in 2-D interpolation using Franke's test function and corresponding optimal values of $\epsilon$, $\alpha$, and $\beta$. The objective function for the optimization was RMS Error. The eigenvalue spectra of the interpolation matrices are shown for the hybrid kernel (top) and the hybrid kernel with polynomial augmentation (bottom). We have used the data for this numerical test from Table \ref{tab:PSO2DFranke}. The eigenvalues for various degrees of freedom, i.e., $N= [ 25, 49, 144, 196, 625, 1296, 2401, 4096]$ have been plotted together. For the hybrid kernel, the spectra has some negative eigenvalues for $N = [25, 49, 144, 196]$. The reason for this is the significantly large values of $\beta$ (see Table 3) representing the dominance of the cubic kernel in the hybridization. We have observed that similar (slightly smaller) RMS errors can also be achieved if we force the values of $\beta$ to be smaller. For other larger degrees of freedoms (excluding the extra data points for polynomial augmentation for the augmented system), i.e., $N=[625,1296,2401,4096]$, all the eigenvalues are positive. The observation here is that if the cubic part is significantly small, the eigenvalues of the proposed hybrid kernel are positive making the kernel positive definite. On the other hand, if we include polynomial terms in the hybrid kernel, the spectrum has some negative part for all degrees of freedom, even with the very small cubic part in the kernel.

\subsection{The objective functions}
\begin{figure}[hbtp]
\centering
\includegraphics[scale=0.5]{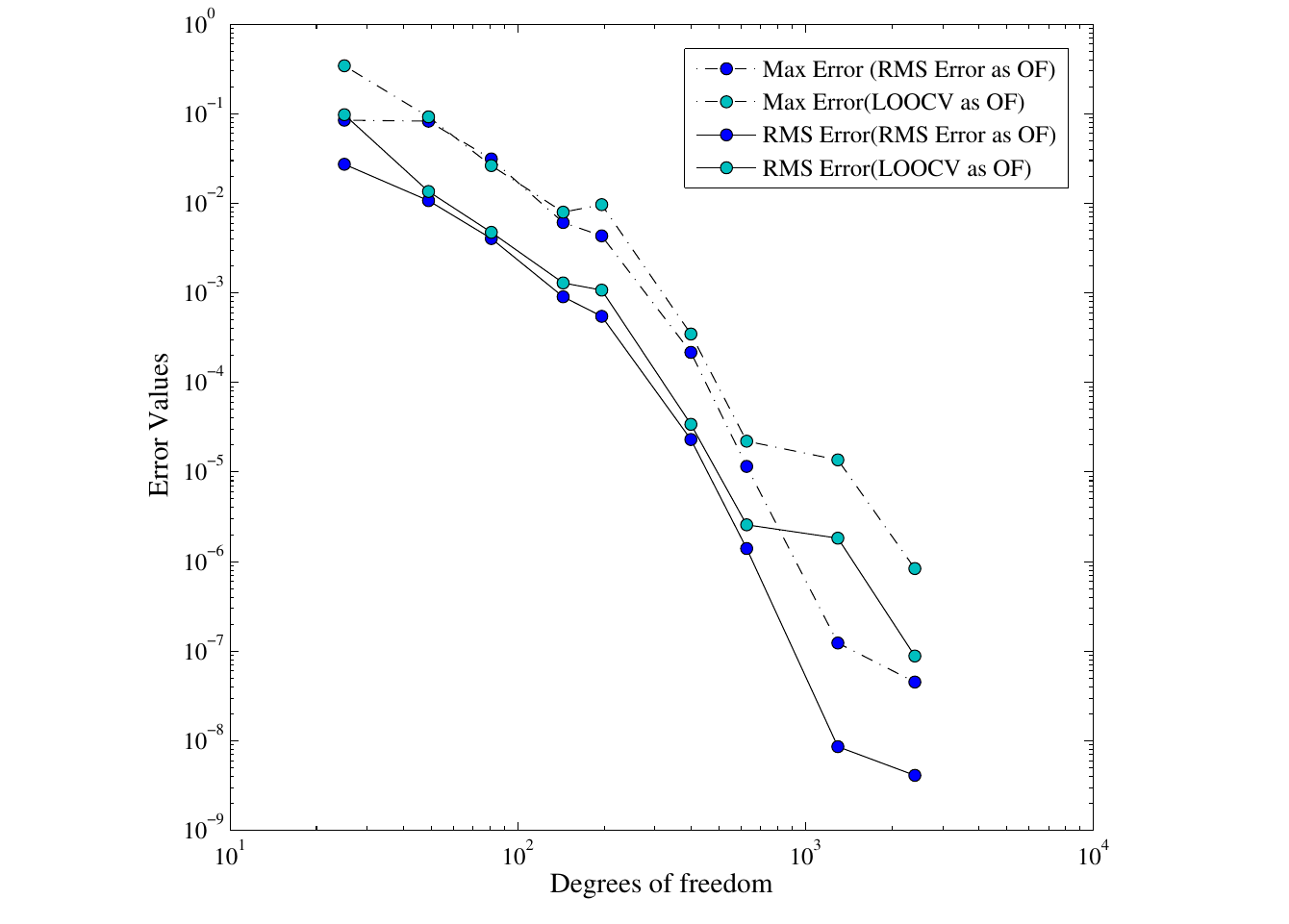}
\caption{Convergence of interpolation errors using two different optimization criteria,\textit{i.e.}, RMS error and LOOCV.}
\label{fig:rmsvsloocv}
\end{figure}

\noindent In this test we compare the optimization for two different objective functions, i.e., the RMS error and the cost function through leave-one-out-cross-validation. We take the same Franke's test function for this interpolation test. The optimization has been performed for various degrees of freedom, i.e., $N= [ 25, 49, 144, 196, 625,\\ 1296, 2401, 4096]$. The results have been tabulated in Table \ref{tab:rmsvsloocv}. $\epsilon$,$\alpha$, $\beta$, and $E_{max}$ are the shape parameter and weight coefficients for the Gaussian and the cubic kernel when the objective function is RMS error, whereas $\epsilon^{'}$,$\alpha^{'}$, $\beta^{'}$ are the similar quantities when the objective function is the cost function given by LOOCV.  Figure \ref{fig:rmsvsloocv} shows the convergence pattern of both the optimizations. The values in Table \ref{fig:rmsvsloocv} suggest that the global minima (as determined by  PSO) are sometimes considerably different, depending on whether one uses RMS error or LOOCV as the objective function. 
\begin{table*}
  \centering \small
  \label{test_kernelp3}
  \begin{tabular*}{\textwidth}{l@{\extracolsep\fill}cccccc}
   \hline
    N & $\epsilon $ & $\alpha $ & $\beta$   & $\epsilon^{'}$ & $\alpha^{'}$ & $\beta^{'}$  \\
    \hline
    25	&2.9432 &0.3161  &0.4660	 &2.4150	 &0.3694	&0.7226	 \\
49	&4.8584	 &0.0770	&0.5830	 &2.7318	 &0.6733	&0.6690	\\
81	&5.1344	 &0.0407	&0.8518	 &3.6403	 &0.7150	&0.06492	\\
144	&6.2932	 &0.0170	&0.8494	 &4.1900	 &0.7419	&0.2815\\
196	&5.2291   &	0.8033	&$3.40e-03$	 &4.3400 &	0.7400	&$2.00e-03$\\
400	&5.5683   &	0.4500	&$4.65e-05$	&5.3350 &0.9203 &$4.56e-08$	\\
625	&5.5434   &	0.6749	&$4.91e-07$	&5.1700 &0.6690 &$1.41e-07$	\\
1296&6.2474 &	0.7880	&$9.11e-09$	&4.2420 &0.9300 &$4.00e-03$	 \\
2401&5.8711 &	0.6318	&$1.28e-07$	&5.1105 &0.7427 &$1.90e-05$	\\

    \hline

    \hline
\end{tabular*}
\caption{The optimized parameters using two different objective functions (OF),\textit{ i.e.}, RMS error and Leave-one-out-cross-validation.}
\label{tab:rmsvsloocv}
\end{table*}

\subsection{Interpolation of normal fault data}
\begin{figure}[hbtp]
\centering
\includegraphics[scale=0.4]{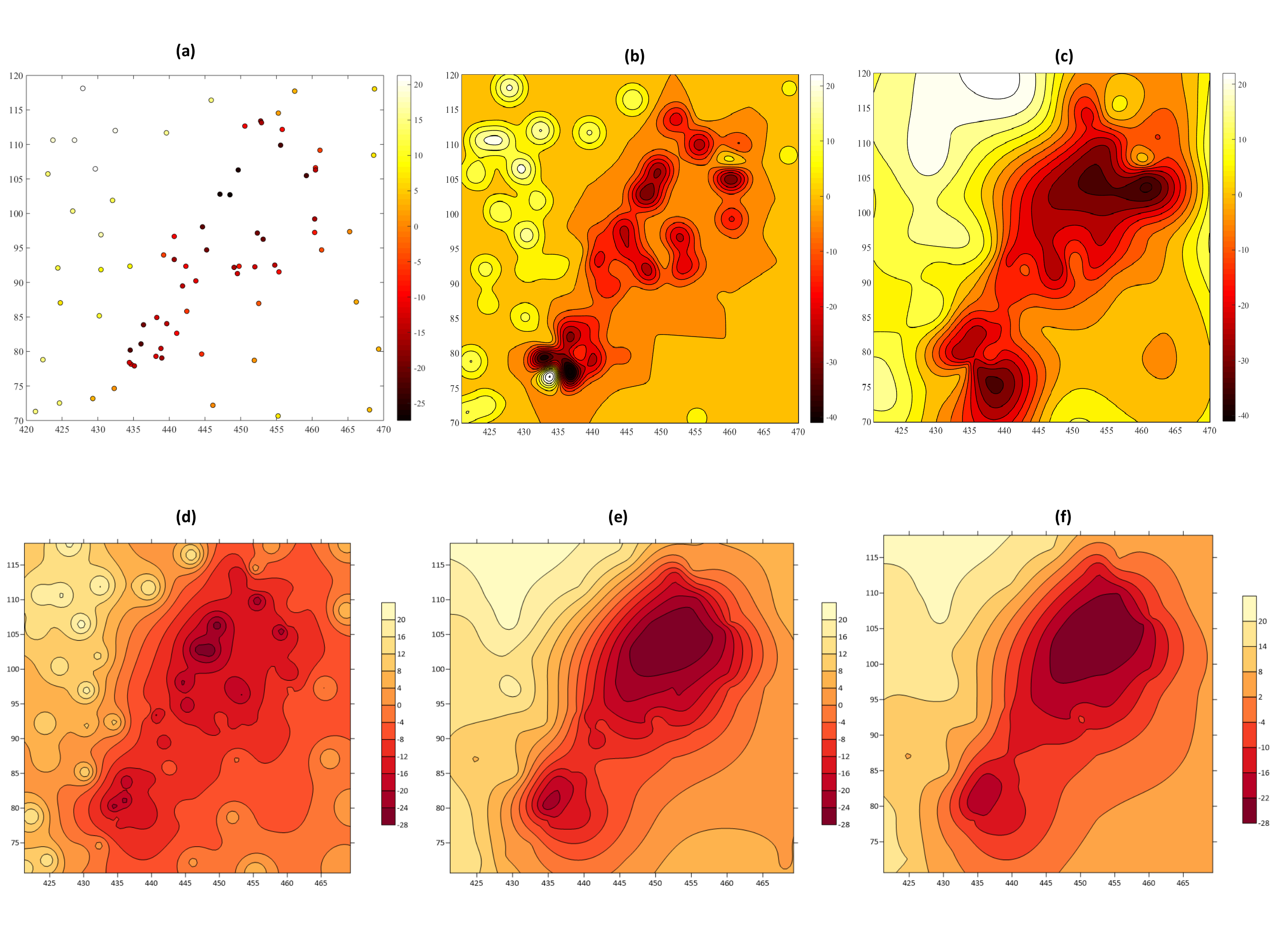}
\caption{(a) The scatter plot of the original data. Interpolation of the normal fault data using various techniques (b) Gaussian RBF (c) Hybrid RBF (d) Inverse Distance (Surfer) (e) Ordinary Kriging (Surfer), and (f) Universal Kriging with linear drift (Surfer).}
\label{fig:compareinterp}
\end{figure}

\begin{figure}
\centering
\includegraphics[scale=0.6]{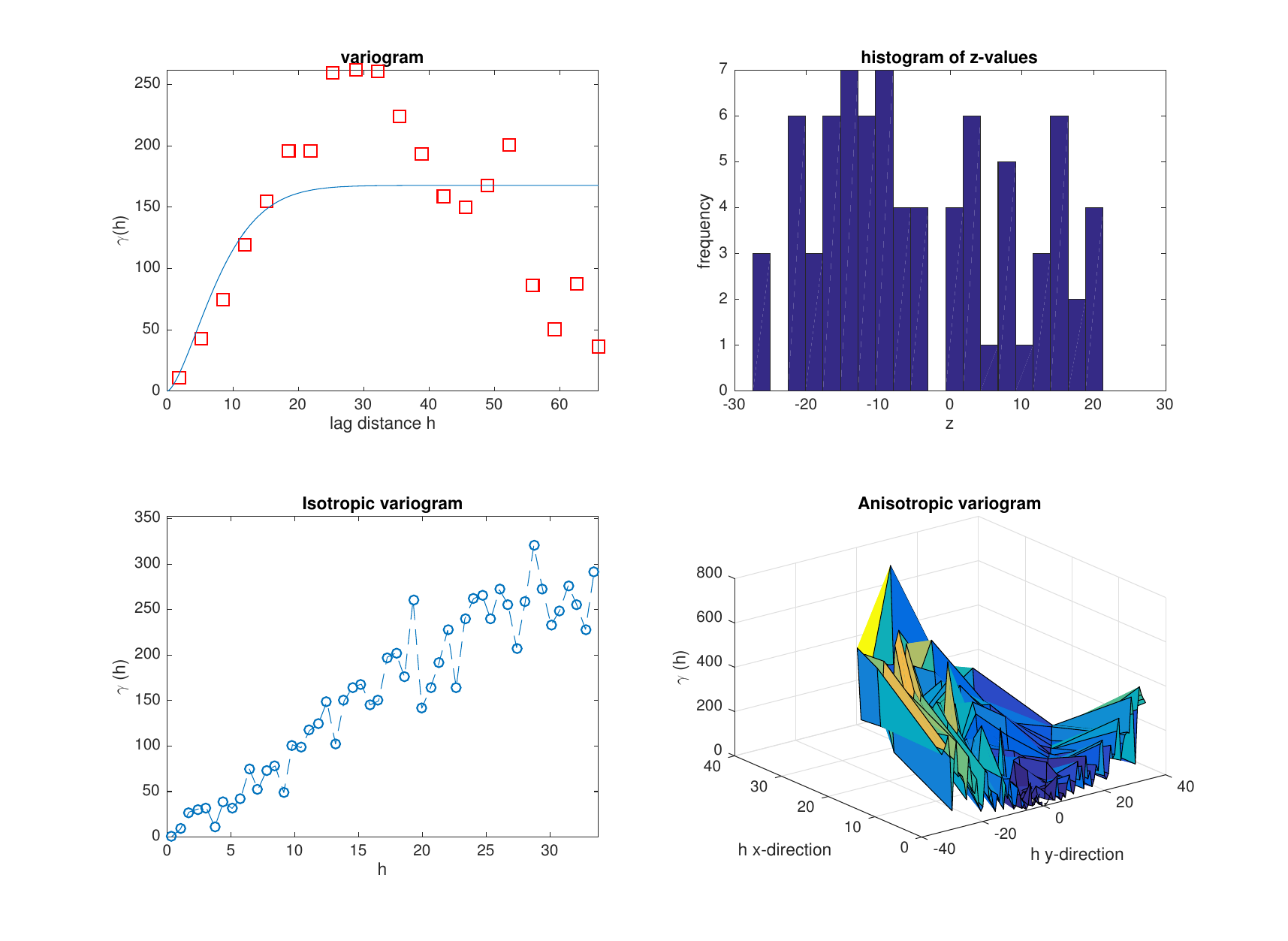}
\caption{Variogram models for the data used in section 6.5. We have used the following MATLAB package for the same https://www.mathworks.com/matlabcentral/fileexchange/29025-ordinary-kriging.}
\label{fig:vario}
\end{figure}

\noindent Next, an interpolation test is performed for synthetic geophysical data. For this, a synthetic dataset representing the vertical distance of the surface of a stratigraphic horizon from a reference surface is considered \citep{Martin2010}. The data set contains a lithological unit displaced by a normal fault. The foot wall has very small variations in elevations whereas the elevation in the hanging wall is significantly variable representing two large sedimentary basins. This data contains the surface information at irregularly spaced $78$ locations in the 2500 $km^2$ domain as shown in Figure \ref{fig:compareinterp}. This data has been reconstructed at $501\times501$, i.e., $251001$ regularly spaced new locations using hybrid radial basis interpolation. Although the used data is synthetic, the exact solution is not known. Figure \ref{fig:vario} shows different variogram model for this data. Optimization has been performed for this case, using the cost function, generated by the LOOCV scheme, which was explained earlier.  The optimal value of the parameters are $\epsilon=0.4318$, $\alpha = 0.7265$, and $\beta=0.4440$. Since we have already observed that the optimum value of $\epsilon$ is similar for the Gaussian kernel and the hybrid kernel, we compare the interpolation outputs with the Gaussian kernel using the same shape parameter as the hybrid one.

Figure \ref{fig:compareinterp} shows the interpolation of this fault data with various interpolation techniques like, MATLAB's linear and cubic interpolation, ordinary and universal kriging, Gaussian RBF and hybrid RBF interpolation. It can be seen that the conventional Gaussian kernel (and Inverse Distance approach) does not interpolate this data properly, and the interpolation using the hybrid kernel provides interpolation in good agreement with the Ordinary and the Universal kriging.

\subsection{Computational Cost}
\begin{figure}[hbtp]
\centering
\includegraphics[scale=0.6]{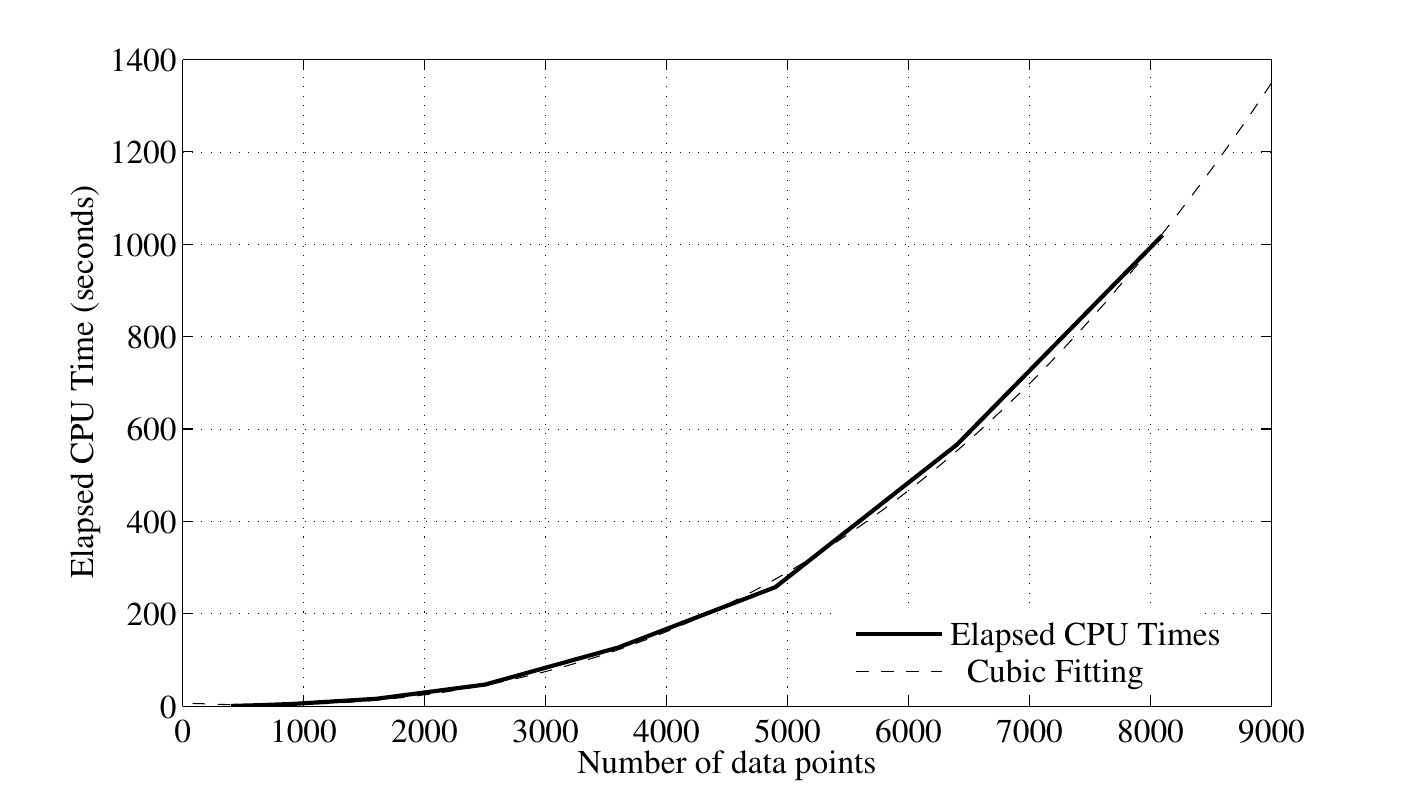}
\caption{Elapsed CPU time 2D RBF interpolation with hybrid Gaussian-cubic kernels, including the time taken in the optimization process. The PSO algorithm in this test used the swarm size of 20 for a single iteration.}
\label{fig:cputime}
\end{figure}

The computational cost study has been performed using a MATLAB implementation of the 2D interpolation algorithm using Franke's test function. Figure \ref{fig:cputime} illustrates the computational cost of 2D interpolation with the proposed hybrid kernel including the cost of parameter optimization. The cost of global RBF interpolation with the hybrid kernels varies approximately as $N^3$, which is similar to the RBF-QR approach without the optimization of kernel parameters \cite{Forn2011}. Moreover, the cost of RBF-QR significantly increases with the optimal value of the shape parameter, which is not the case here. The computational cost of the presented approach is likely to be reduced when used in local approximation form such as radial basis-finite difference (RBF-FD)--- as the ``matrices in the RBF-FD methodology go from being completely full to  $99\%$ empty" \citep{Flyer2020}.

\section{Conclusions}
\noindent We have proposed a hybrid kernel for radial basis interpolation and its applications. This hybrid kernel utilizes the advantages of the Gaussian and the cubic kernel according to the problem type. Based on the numerical tests performed in this work, we draw the following conclusions:
\begin{enumerate}
\item Combination of a small part of the cubic kernel in the Gaussian kernel reduces the condition number significantly, making the involved algorithm well-posed.
\item The optimal value of the shape parameter remains nearly the same for the Gaussian and hybrid Gaussian-cubic kernel.
\item The interpolation using the proposed hybrid kernel remains stable under the low shape parameter paradigm unlike the one with only the Gaussian RBF.
\item The hybrid kernel was used to interpolate real type geophysical data, which had close observational points and large changes, due to which, the Gaussian kernel could not interpolate this data. However, the interpolation using the hybrid kernel exhibits convincing results in agreement with the ordinary and universal kriging approach.
\item When used in the global form, the computational cost of the proposed approach was found to vary as $N^3$, which is similar to the RBF-QR approach. However, unlike RBF-QR, the cost of the present approach is not parameter dependent. Also, the computational cost of the present approach can be further minimized by normalizing the hybrid kernel---therefore reducing the  number of kernel parameters: from three to two.
\item Future work shall involve the application of the proposed hybrid kernels in local RBF interpolations such as RBF-FD and for stable meshless schemes for numerical solution of PDEs \cite{mishra2017improved,mishra2017frequency}.
\item In this paper we have focused on an interpolation method for (exact) data. If there is noise present in the data, then RBF interpolation methods can be regularized in a manner that is analogous to smoothing splines (or ridge regression).
\end{enumerate}

\section*{Appendix A: LOOCV}
\noindent Following the notations used in problem 2.1, let us write the datasites without the $k^{th}$ data as,
\[\bm{x}^{[k]} = [\bm{x}_1,...,\bm{x}_{k-1}, \bm{x}_{k+1},...,\bm{x}_N ]^T.\]
The removed point has been indicated by the superscript $[k]$. This superscript will differentiate the quantities computed with ``full" dataset and partial data set without the $k^{th}$ point. Hence, the partial RBF interpolant $\mathcal{F}(\bm{x})$ of the given data $\bm{f}(\bm{x})$ can be written as,
\[\mathcal{F}(\bm{x}) = \sum_{j=1}^{N-1} c^{[k]}_{j} \phi (\parallel \bm{x}-\bm{x}^{[k]}_{j}\parallel). \]
The error estimator can, therefore, be written as,
\[e_k = \bm{f}(\bm{x}_k) - \mathcal{F}^{[k]}(\bm{x}_k).\]
The norm of the error vector $\bm{e} = [ e_1,..., e_N]^T$, obtained by removing each one point and comparing the interpolant to the known value at the excluded point determines the quality of the interpolation. This norm serves as the ``cost function" which is the function of the kernel parameters $\epsilon$, $\alpha$, and $\beta$. We consider $l_2$ norm of the error vectors for our purpose. The algorithm for constructing the ``cost function" for RBF interpolation via LOOCV has been summarized in Algorithm 1. We recommend \cite{Fass2009,Fass2012} for some more insights of the application of LOOCV in radial basis interpolation problems. Here $c_k$ is the $k^{th}$ coefficient for the interpolant on ``full data" set and $\mathbf{A}^{-1}_{kk}$ is the $k^{th}$ diagonal element in the inverse of the interpolation matrix for ``full data".
\begin{algorithm}[!htbp]
\begin{algorithmic}[1]
\STATE Fix a set of parameters $[\epsilon, \alpha, \beta]$
\STATE For $k=1,...,N$
\STATE Compute the interpolant by excluding the $k^{th}$ point as,
\begin{eqnarray}
\mathcal{F}(\bm{x}) = \sum_{j=1}^{N-1} c^{[k]}_{j} \phi (\parallel \bm{x}-\bm{x}^{[k]}_{j}\parallel).
\end{eqnarray}
\STATE Computer the $k^{th}$ element of the error vector $e_k$
\begin{eqnarray}
e_k = \mid\bm{f}(\bm{x}_k) - \mathcal{F}^{[k]}(\bm{x}_k) \mid,
\end{eqnarray}
\STATE As proposed by Rippa \citep{Rippa}, a simplified alternative approach to compute $e_k$ is,
\begin{eqnarray}
e_k =  \frac{c_k}{\mathbf{A}^{-1}_{kk}}.
\end{eqnarray}
\STATE end
\STATE Assemble the ``cost vector'' as $\bm{e} = [ e_1,..., e_N]^T$.
\end{algorithmic}
\caption{Leave-one-out-cross-validation for radial basis interpolation schemes.}
\label{alg:test}
\end{algorithm}

\section*{Appendix B: Particle swarm optimization}
\label{AA}
\begin{figure}[hbtp]
\centering
\includegraphics[scale=0.6]{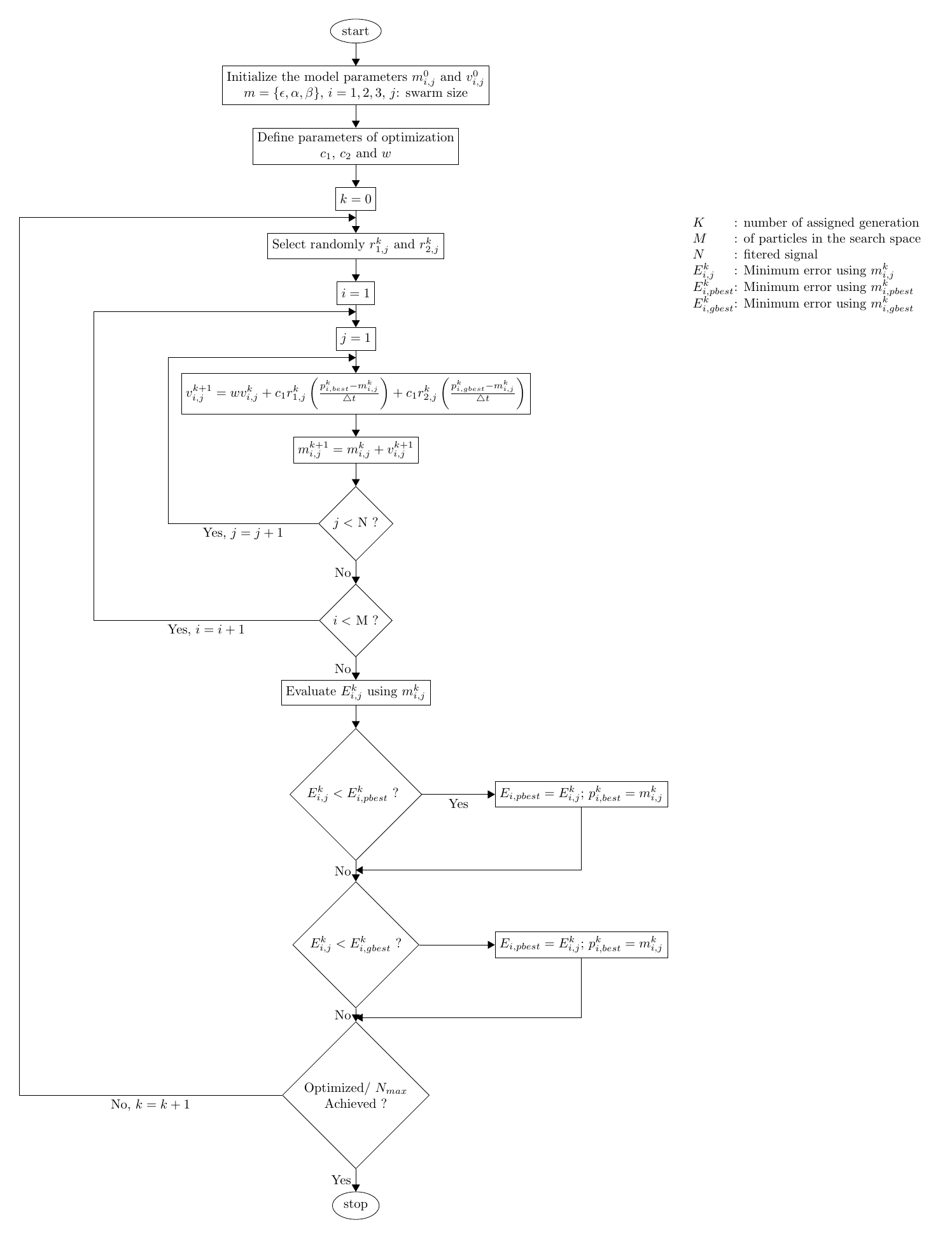}
\caption{Flowchart of particle swarm optimization in the context of numerical tests.}
\label{fig:psoflowchart}
\end{figure}

\noindent The term optimization refers to the process of finding a set of parameters corresponding to a given criterion among many possible sets of parameters. One such optimization algorithm is particle swarm optimization (PSO), proposed by James Kennedy and Russell Eberhart in 1995 \citep{Everhart1995,Everhart2001}. PSO is known as an algorithm which is inspired by the exercise of living organisms like bird flocking and fish schooling. In PSO, the system is initiated with many possible random solutions and it finds optima in the given search space by updating the solutions over the specified number of generations. The possible solutions corresponding to a user defined criterion are termed as \emph{particles}. At each generation, the algorithm decides optimum particle towards which all the particles fly in the problem space. The rate of change in the position of a particle in the problem space is termed as \emph{particle velocity}. In each generation, all the particles are given two variables which are known as \emph{pbest} and \emph{gbest}. The first variable (\emph{pbest}) stores the best solution by a particle after a typical number of iteration. The second variable (\emph{gbest}) stores the global best solution, obtained so far by any particle in the search space \cite{Shaw2007,Singh2015}. Once the algorithm finds these two parameters, it updates the velocity and the position of all the particles according to the following pseudo-codes,
\begin{eqnarray}
\nonumber
v[.]= v[.]+c_1*rand(.)*(pbest[.] - present[.]) + c_2 * rand(.) * (gbest[.] - present[.]),
\end{eqnarray}
\begin{eqnarray}
\nonumber
present[.] = present[.] + v[.]
\end{eqnarray}
Where, $v[.]$ is the particle velocity, present[.] is the particle at current generation, and $c_1$ and $c_2$ are learning factors. According to the studies of Perez and Behdinan \cite{Perez2007}, the particle swarm algorithm is stable only if the following conditions are fulfilled;
\[0< c_1+c_2 < 4\]
\[\left( \frac{c_1+c_2}{2}\right)-1 < w<1\]
The optimization of the parameters of hybrid Gaussian-cubic kernel using particle swarm optimization is summarized in the flowchart given in the Figure \ref{fig:psoflowchart}.

\section*{Appendix C: Variogram models of `normal fault' data}

\bibliographystyle{spmpsci}       
\bibliography{sample}   
\end{document}